\documentclass[oneside,11pt,reqno]{amsart}
\usepackage{amssymb,amsmath,amsthm,bbm,enumerate,mdwlist,url,multirow,hyperref,amsthm,mathabx}
\usepackage[shortlabels]{enumitem}

\theoremstyle{definition}
\newtheorem{definition}{Definition}%Extra square-bracket argument achives that the numbering is the same as for definition (single uniform counter). 
\theoremstyle{theorem}
\newtheorem{proposition}[definition]{Proposition}
\newtheorem{theorem}[definition]{Theorem}
\theoremstyle{remark}
\newtheorem{remark}[definition]{Remark}

\def\PP{\mathsf P}

\def\EE{\mathsf E}

\def\LL{\mathcal L}
\def\FF{\mathcal F}

\def\Ber{\mathrm{Ber}}
\def\Exp{\mathrm{Exp}}
\def\geom{\mathrm{geom}}
\def\Pois{\mathrm{Pois}}
\def\HPP{\mathrm{HPP}}
\def\fii{\varphi}
\def\cov{\mathrm{cov}}
\begin{document}
\title{Another characterization of homogeneous Poisson processes}

\author{Matija Vidmar}

\address{Department of Mathematics, University of Ljubljana, Slovenia}
\address{Institute for Mathematics, Physics and Mechanics, Ljubljana, Slovenia}

\thanks{The author acknowledges financial support from the Slovenian Research Agency (research core funding No. P1-0222). I thank Janez Bernik for suggesting me the investigation of this problem.}
\email{matija.vidmar@fmf.uni-lj.si}


\begin{abstract}
For a general renewal process $N$ (allowing delay, defect and multiple simultaneous arrivals) the independence of the first renewal epochs of the marked processes got from $N$ by Bernoulli $0$/$1$ thinning is characterized. This independence is well-known to hold true in the case of homogeneous Poisson processes; by way of corollary one obtains the interesting observation that, when coupled with some minimal extra conditions, it in fact already identifies them. %a related characterization of homogeneous Poisson processes (HPPs), which answers the natural question of whether or not the described independence property.
The proof is analytic in character.
\end{abstract}


\maketitle

\section{Introduction and problem delineation}\label{section:intro}
A process $N=(N_t)_{t\in [0,\infty)}$ in continuous time is a homogeneous Poisson process (HPP) of intensity $c\in (0,\infty)$, by definition, if it is a counting process (i.e. if it has values in $\mathbb{N}_0\cup\{\infty\}$ and right-continuous nondecreasing paths) that is finite a.s., has jumps of size $1$ a.s., starts in $N_0=0$ a.s., and has independent increments, whose distribution is Poisson: $(N_t-N_s)\mathbbm{1}_{\{N_s<\infty\}}\sim \Pois(c(t-s))$ for $ \{s,t\}\subset [0,\infty)$, $s\leq t$. (Here, for $\lambda\in [0,\infty)$, $\Pois(\lambda)$ is the law on $\mathbb{N}_0$ that has $\Pois(\lambda)(\{k\})=\lambda^ke^{-\lambda}/k!$ for $k\in \mathbb{N}_0$.) We write: $N\sim \HPP(c)$. Homogeneous Poisson processes represent a fundamental type (in law a one-parametric family) of processes in continuous time and with a discrete state space, lying on the intersection of (at least) counting processes, L\'evy processes (hence strong Markov processes), renewal processes and (inhomogeneous) Poisson processes (hence continuous-time Markov chains). Consequently there has been, and there exists, a considerable interest in various characterizations of HPPs. An incomplete but illustrative list of such characterizations follows. (For unexplained terms the reader is referred to the cited works.)

For $c\in (0,\infty)$ and a counting process $N$ with $N_0=0$ a.s.,  the assertion ``$N\sim \HPP(c)$'' is equivalent to each of the following:
\begin{itemize}
\item $N$ is a  L\'evy process with jumps of size $1$ a.s. and with its L\'evy measure having mass $c$ \cite[Theorem~2.2.13]{applebaum}.
\item $N$ has jumps of size $1$ a.s. and its inter-arrival times are independent, identically, exponentially with mean $c^{-1}$  (notation: $\Exp(c)$), distributed \cite[Theorem~6.5.5(d)]{cinlar}.
\item (Watanabe) $N$ has jumps of size $1$ a.s. and $(N_t-c t)_{t\in [0,\infty)}$ is a martingale \cite[Theorem~6.5.5(c)]{cinlar}.
\item $N$ is an ordinary  (i.e. non-delayed non-defective no-simultaneous-arrivals) renewal process that is stationary, with its inter-renewal times having mean $c^{-1}$: follows from \cite[Corollary~V.3.6]{asmussen} coupled with the elementary observation that invariance under the transformation of the integrated tail characterizes the exponential distribution.
%\item $N$ is an (inhomogeneous) Poisson process with intensity function $[0,\infty)\in t\mapsto ct\in [0,\infty)$. \cite[Theorem~2.1.1 (proof thereof)]{ross}
\item (Srivastava) A non-trivial Bernoulli marking (thinning) of $N$ results in independent marked  processes and $\EE N_t=ct$ for $t\in [0,\infty)$: this is a particular case of the more general characterization of Poisson processes in Euclidean space \cite[Theorem~2.1]{assuncao}, originally due to Fichtner; see \cite[Theorem~1]{nehring} for further extensions.
\item (Samuels) $N$ is an ordinary renewal process with mean inter-renewal time $c^{-1}$ that results as the superposition of two independent ordinary renewal processes \cite[Theorem on p. 73]{samuels}. \label{samuels}
\end{itemize}
For still further characterizations see \cite[Sections~2.2,~2.3 and~4.3]{vere-jones} dealing with the property of `complete randomness', distribution form and various operations on stationary renewal processes, respectively; Slivnyak-Mecke's characterization of Poisson processes \cite[Proposition 13.1.VII]{vere-jones-II} \cite[Lemma~6.15]{kallenberg} in the context of Palm theory of random measures; \cite{liberman,gan} for the order statistics property; \cite{jagers,li,erickson} concerning age (a.k.a. spent or current life) and residual life; finally  \cite{huang,chandramohan,disney} that deal with marking (thinning) of renewal processes. \label{lit}

In this paper we present another characterization of HPPs in the context of marking  (thinning) a general renewal process.  To this end we first fix some notation. 

Let, on a probability space $(\Omega,\FF,\PP)$, $T=(T_i)_{i\in \mathbb{N}}$ be a sequence of independent random variables with values in $[0,\infty]$. Let $T_j$, $j\in \mathbb{N}_{\geq 2}$, be identically distributed. Define $S_n:=\sum_{i=1}^nT_i$ for $n\in \mathbb{N}$, and then $N_t:=\sum_{n\in \mathbb{N}}\mathbbm{1}(S_n\leq t)$ for $t\in [0,\infty)$ -- the associated renewal process (allowing delay: $T_1$ does not necessarily have the same distribution as $T_2$; defect: $T_i$, $i\in \mathbb{N}$, can take on the value $\infty$; and multiple simultaneous arrivals: $T_i$, $i\in \mathbb{N}$, can take on the value $0$).

Let furthermore $p\in (0,1)$, and let $X=(X_i)_{i\in \mathbb{N}}$ be a sequence of independent, identically distributed (i.i.d.) random variables taking values in $\{0,1\}$, independent of  $T$, and with $X_1\sim \Ber(p)$, where $\Ber(p)$ is the Bernoulli law: $\Ber(p)(\{1\})=1-\Ber(p)(\{0\})=p$. Define the marked processes $N^1$ and $N^0$ as follows: $$N_t^i:=\sum_{n\in \mathbb{N}}\mathbbm{1}(S_n\leq t,X_n=i),\quad t\in [0,\infty),\quad i\in \{0,1\}.$$ The strong Markov property for i.i.d. sequences implies that $N^0$ and $N^1$ are again renewal processes (with delay and defect): for $i\in \{0,1\}$, if one defines $S^i_0:=0$ and inductively $S^i_{n+1}:=\inf\{m>S^i_n:X_m=i\}$, $n\in \mathbb{N}_0$, then $(S^i_{n+1}-S^i_n)_{n\in \mathbb{N}_0}$ (where we set e.g. $\infty-\infty=\infty$ on the negligible event on which such a difference may occur) is an i.i.d. sequence, independent of $T$, and a.s. the sequence of the inter-renewal epochs of $N^i$ is given by $(\sum_{k=S^i_{n-1}+1}^{S^i_n}T_k)_{n\in \mathbb{N}}$. 

Finally define, for $i\in \{0,1\}$, $R_i:=\inf\{t\in [0,\infty):N^i_t\geq 1\}$, the time of the first renewal of the process $N^i$, and $L_i:=\inf\{j\in \mathbb{N}:X_j=i\}$. Then, on $\{L_i<\infty\}$ and hence a.s., $R_i=\sum_{j=1}^{L_i}T_j$, $i\in \{0,1\}$.\label{ns}

We observe: if $N\sim \HPP(\theta)$ for some $\theta\in (0,\infty)$, then $N^0$ and $N^1$, in particular $R_0$ and $R_1$, are independent \cite[Theorem~4.4.1]{resnick}. It is then natural and interesting to ask, whether or not the latter property of the independence of $R_0$ and $R_1$ already characterizes HPPs. Indeed we will demonstrate the validity of\label{natural}

\begin{theorem}\label{theorem}
Assume $\PP(T_1<\epsilon)>0$ for all $\epsilon>0$, and that either $T_2$ is non-arithmetic or else  $\PP(T_1=0)=0$. Then  $R_0$ and $R_1$ are independent if and only if $N\sim\HPP(\theta)$ for some $\theta\in (0,\infty)$. The same equivalence obtains if $N$ is assumed to be ordinary (i.e. non-delayed, non-defective, and not having multiple simultaneous arrivals) instead. 
\end{theorem}
Here:
\begin{definition}\label{definition}
$T_2$ is non-arithmetic, if there is no $\alpha\in (0,\infty)$ with $\PP(T_2\in \{\alpha n:n\in \mathbb{N}_0\cup \{\infty\}\})=1$.
%\leavevmode
%\begin{enumerate}
%\item Nearitmetičnost $T_2$ pomeni, po definiciji, da ne obstaja $\alpha\in (0,\infty)$ z $\PP(T_2\in \{\alpha n:n\in \mathbb{N}_0\cup \{\infty\}\})=1$.
%\item Modulo trivialen determinističen zamik $\kappa$, je to natančna karakterizacija HPP. 
%\end{enumerate}
\end{definition}
Theorem~\ref{theorem}, whose proof is given at the end of Section~\ref{section:result}, is most closely related to the findings of \cite{chandramohan,disney,huang}. Let us see how it compares.  On the one hand, \cite[Corollaries~2.3 and~3.2]{disney} (respectively, \cite[Theorem~2.1]{chandramohan};  \cite[Corollary~2]{huang}) give that, when either $T_1$ has the distribution of the integrated tail of $T_2$ with (implicitly) $\EE T_2\in (0,\infty)$, or else when $T_1$ has the same distribution as $T_2$, $\PP(T_2<\infty)=1$  and $T_2$ is non-arithmetic (respectively, when $\PP(0<T_1)=1$ and (as implicitly assumed in the proof; not all the assumptions appear to be given explicitly) $\PP(T_1=\infty)<1$; when $\PP(0<T_1,0<T_2)=1$ and $\PP(T_1<\epsilon)>0$ for all $\epsilon>0$), then $\cov(N^0_t,N^1_t)=0$ for all $t\in (0,\infty)$ implies $N\sim \HPP(\theta)$ for some $\theta\in (0,\infty)$.  (Strictly speaking the quoted result of \cite{chandramohan} is false. For, given a $\kappa\in (0,\infty)$, we can take independent $T_1-\kappa\sim \Exp(\lambda)$ and $T_j\sim \Exp(\lambda)$ for $j\in \mathbb{N}_{\geq 2}$ (a deterministically delayed HPP). This situation is however precluded by the conditions of \cite{disney,huang}.) On the other hand, the condition of Theorem~\ref{theorem} is one on the independence of \emph{the first renewal epochs $R^0$ and $R^1$ only}, and \emph{not} (\emph{a priori}) on the absence of correlation of the processes $N^0$ and $N^1$ at all deterministic times (viz. the condition of  \cite{chandramohan,disney,huang}). In a similar vein, Theorem~\ref{theorem} is not subsumed in the result of Samuels described above (final bullet point on p.~\pageref{samuels}): in the case that $N$ is an ordinary renewal process, for sure $N^0$ and $N^1$ are ordinary renewal processes that superpose into $N$, but the condition of Theorem~\ref{theorem} is \emph{not} (\emph{a priori}) on $N^0$ and $N^1$ being independent. Thus Theorem~\ref{theorem} is a complement to existing characterizations of HPPs in the context of marked renewal processes. 

In fact we shall prove slightly more than what is the contents of Theorem~\ref{theorem}. Specifically, we shall provide a precise characterization of the independence of $R_0$ and $R_1$ (see Proposition~\ref{proposition:karakterizacija} in the section following),  whose immediate corollary will be Theorem~\ref{theorem}. Excepting degenerate and trivial cases, and modulo  deterministic time delay and scaling, we obtain here besides HPPs also what are continuous-time-embedded discrete-time stationary  ``geometric'' renewal processes. It is interesting that these yield independence of the first renewal epochs of the two marked processes, however not the independence of the marked processes in their entirety (see Remark~\ref{remark}\ref{remark:c}).

In addition to its theoretical appeal, our result appears to have some potential practical (statistical) relevance as well. We mean here a situation in which, for some reason, $N$ may be assumed to satisfy the assumptions of Theorem~\ref{theorem}, but it is not clear whether $N$ is an HPP. Then this can be statistically tested, via independent trials, based on (also) the hypothesis of the independence of the first renewal times $R_0$ and $R_1$ of a non-trivial Bernoulli marking (possibly of unknown parameter $p$) of $N$. %Moreover, in the non-defective case, a statistics on this can be gathered, by independent trials, in a finite (though random) time. 
This might in particular be useful when data is limited to $R_0$ and $R_1$. Compare the study \cite{assuncao} of a test for Poisson processes, based on the characterization result of Fichtner alluded to above. On a more pure level, note that $R_0$ and $R_1$ are relatively simple functionals of the paths of $N^0$ and $N^1$, and as such, in some given context, their independence may be more easily susceptible to analysis, than that of the whole of the processes $N^0$ and $N^1$, or of $N$. \label{practical}

\section{The result and its proof}\label{section:result}
So as to be able to state the precise result of this paper succinctly, let us agree on the following pieces of notation: $\LL(V)$ denotes the law of a random element $V$; for $x_0\in [0,\infty]$, $\delta_{x_0}$ is the Dirac measure at $x_0$; then for $r\in (0,1)$,  $\geom_{\mathbb{N}}(r):=\sum_{k=1}^\infty r(1-r)^{k-1}\delta_k$, respectively $\geom_{\mathbb{N}_0}(r):=\sum_{k=0}^\infty r(1-r)^{k}\delta_k$, is the geometric law on $\mathbb{N}$, respectively $\mathbb{N}_0$, with success parameter $r$.

Now the result of this paper follows. 

\begin{proposition}\label{proposition:karakterizacija}
$R_1$ and $R_0$ are independent if and only if (precisely) one of the conditions below holds true. 
\begin{enumerate}[(a)]
\item\label{kara:4}  $\PP(T_1=\infty)=1$.
\item\label{kara:3} There exists  $\kappa\in [0,\infty)$ such that $\PP(T_2=0)=\PP(T_1=\kappa)=1$.
\item\label{kara:0}  There exist $\kappa\in [0,\infty)$ and $q_0\in (0,1)$ such that $\LL(T_1)=(1-q_0^2)\delta_\kappa+q_0^2\delta_\infty$ and $\LL(T_2)=(1-q_0)\delta_0+q_0\delta_\infty$.%$\PP(T_2=0)=\PP(T_1=\kappa)=1$;
\item\label{kara:2} There exist $\kappa\in [0,\infty)$ and $\theta\in (0,\infty)$ such that $T_1-\kappa\sim \Exp(\theta)$ and $T_2\sim \Exp(\theta)$.
\item\label{kara:1} There exist $q_0\in (0,1)$, $\kappa\in [0,\infty)$ and $\alpha\in (0,\infty)$ such that $\LL((T_1-\kappa)/\alpha)=\geom_{\mathbb{N}_0}(1-q_0^2)$ and $\LL(T_2/\alpha)=(1-q_0)\delta_0+q_0\geom_\mathbb{N}(1-q_0^2)$.
\end{enumerate}
\end{proposition}

\begin{remark} \label{remark}
\leavevmode
\begin{enumerate}[(i)] 
\item The conditions of the proposition are clearly mutually exclusive. 
\item\label{remark:ii} In cases \ref{kara:4}, \ref{kara:3} and \ref{kara:2} even the processes $N^1$ and $N^0$ in their entirety are independent.
\item\label{remark:c} In cases \ref{kara:1} and \ref{kara:0}, $N^1$ and $N^0$ are not independent. This may be seen as follows. %Označimo z $L^i$ markirana procesa diskretnega prenovitvenega procesa $L$: $$ L_n^i:=\sum_{k=1}^\infty\mathbbm{1}(\sum_{l=1}^kU_l\leq n,X_k=i)\quad n\in \mathbb{N}_0,\quad i\in \{0,1\}.$$ 
Let $B_i:=\{N^i_\kappa=1\}$, $i\in \{0,1\}$. We compute $\PP(B_0)=\sum_{k=1}^\infty (1-q_0^2)(1-q_0)^{k-1}q_0{k\choose 1}p^{k-1}(1-p)=q_0(1-q_0^2)(1-p)/(1-(1-q_0)p)^2$, and similarly $\PP(B_1)=q_0(1-q_0^2)p/(1-(1-q_0)(1-p))^2$, finally $\PP(B_1\cap B_0)=2(1-q_0^2)(1-q_0)q_0p(1-p)$. Let furthermore $A_0:=\{N^0_\kappa=0\}$. We compute $\PP(A_0)=q_0^2+\sum_{k=1}^\infty (1-q_0^2)(1-q_0)^{k-1}q_0p^k=q_0^2+(1-q_0^2)q_0p/(1-(1-q_0)p)=q_0(q_0+p-pq_0)/(1-p+pq_0)$ and also $\PP(A_0\cap B_1)=(1-q_0^2)pq_0$. 
Elementary simplifications reveal that $\PP(B_0\cap B_1)=\PP(B_0)\PP(B_1)$ is equivalent to $q_0(1+q_0)=2[q_0+p-pq_0]^2[1-p+q_0p]^2$,  whilst $\PP(A_0\cap B_1)=\PP(A_0)\PP(B_1)$ is equivalent to $[q_0+p-pq_0][1-p+q_0p]=q_0$. Both equalities together yield $q_0(1+q_0)=2q_0^2$, contradicting $q_0\in (0,1)$.%, in nato prva $q_0^4(1-q_0^2)=[q_0+p-pq_0]^2[1-p+q_0p]^2$, vendar pa imamo v slednji namesto $=$, $<$, ker je $q_0^2(1-q_0^2)<[1-p+q_0p]^2$, saj desna strani slednje neenakosti doseže svoj minimum pri $p=1$. 

% Naj bo naprej $A_i:=\{N^i_\kappa=0\}$, $i\in \{0,1\}$. Izračunamo $\PP(A_0)=q_0^2+\sum_{k=1}^\infty (1-q_0^2)(1-q_0)^{k-1}q_0p^k=q_0^2+(1-q_0^2)q_0p/(1-(1-q_0)p)$ 
%, in podobno 
%$\PP(A_1)=q_0^2+(1-q_0^2)q_0(1-p)/(1-(1-q_0)(1-p))$ %, končno $\PP(A_0\cap A_1)=q_0^2$. % 

 %Iz zadnje enakosti sledi $q_0=1/\sqrt{2}$, kar da protislovje v prvi enakosti (funkcija $[0,1]\ni r \mapsto (1/\sqrt{2})^3-2[((1/\sqrt{2})+r-r/\sqrt{2})(1-r+r/\sqrt{2})]^2(1-1/\sqrt{2})$ doseže svoj globalni minimum v $r=1/2$ in je tam strogo pozitivna, kot je trivialno preveriti). 

%Elementarno preurejanje razkrije, da je $\PP(B_0\cap B_1)=\PP(B_0)\PP(B_1)$ ekvivalentno $q_0^3=2[q_0+p-pq_0]^2[1-p+q_0p]^2(1-q_0)$,  ter da je $\PP(A_0\cap A_1)=\PP(A_0)\PP(A_1)$ ekvivalentno $(-1 + p) p q_0^2 (-1 + 2 q_0^2)=0$. Iz zadnje enakosti sledi $q_0=1/\sqrt{2}$, kar da protislovje v prvi enakosti (funkcija $[0,1]\ni r \mapsto (1/\sqrt{2})^3-2[((1/\sqrt{2})+r-r/\sqrt{2})(1-r+r/\sqrt{2})]^2(1-1/\sqrt{2})$ doseže svoj globalni minimum v $r=1/2$ in je tam strogo pozitivna, kot je trivialno preveriti). 
%\item Also in case \ref{kara:0} the processes $N^1$ and $N^0$ are not independent. This is seen in precisely the same manner as in the previous item (the calculations are identical!). 
\item The last case (item ~\ref{kara:1}) is (modulo the deterministic time delay by $\kappa$ and the scaling $\alpha$) a stationary renewal process in discrete time that has been embedded into continuous time. For, if $U_1\sim \geom_{\mathbb{N}_0}(1-q_0^2)$ and $U_k\sim (1-q_0)\delta_0+q_0\geom_\mathbb{N}(1-q_0^2)$ for $k\in \mathbb{N}_{\geq 2}$ are independent, then we can easily convince ourselves that, for $k\in \mathbb{N}_0$, $\PP(U_1=k)=\PP(U_2>k)/\EE U_2$ ($U_1$ has the distribution of the  ``summed tail'' of $U_2$), and consequently (e.g. via moment functions) that the expected number of renewals at time $k\in \mathbb{N}_0$, i.e.  $\EE \sum_{n=1}^\infty\mathbbm{1}(\sum_{l=1}^nU_l=k)$, is equal to $1/\EE U_2=(1-q_0^2)/q_0$, and hence does not depend on $k$ (property of stationarity; cf. \cite[pp. 34--36]{limnios}). For completeness' sake: the renewal process $L$ in discrete time, associated to the sequence $(U_k)_{k\in \mathbb{N}}$, is of course given by $L_n:=\sum_{k=1}^\infty\mathbbm{1}(\sum_{l=1}^kU_l\leq n)$ for $n\in \mathbb{N}_0$. 
\end{enumerate}
\end{remark}
The proof of Proposition~\ref{proposition:karakterizacija} will be via Laplace transforms. Let us recall this notion. 
\begin{definition}
For a law $\LL$ on the Borel subsets of $[0,\infty]$ we define the Laplace transform of $\LL$ as the function $[0,\infty)\ni \lambda\mapsto \int_{[0,\infty)}e^{-\lambda x}\LL(dx)\in [0,\LL([0,\infty))]$. We denote it by $L_\LL$. 
\end{definition}
\begin{remark}
Such a Laplace transform is by bounded convergence continuous with limit $\lim_\infty L_\LL=\LL(\{0\})$ at $\infty$, its value at $0$ is $L_\LL(0)=\LL([0,\infty))$, and it is nonincreasing. It also determines the law (``injectivity of the Laplace transform''): if for two laws $\LL_1$ and $\LL_2$ on the Borel subsets of $[0,\infty]$, $L_{\LL_1}\vert_{[a,\infty)}=L_{\LL_2}\vert_{[a,\infty)}$ for some $a\in [0,\infty)$, then $\LL_1=\LL_2$.  %Če $\LL_1$ in $\LL_2$ zamenjamo z $e^{-a\cdot}\mathbbm{1}_{[0,\infty)}\cdot \LL_1$ in $e^{-a\cdot}\mathbbm{1}_{[0,\infty)}\cdot \LL_2$, lahko predpostavimo brez škode za splošnost da je $a=0$.  
For, the restrictions of the measures $\LL_1(\cdot\cap [0,\infty))$ and $\LL_2(\cdot\cap [0,\infty))$ to the Borel subsets of $[0,\infty)$ then have the same (finite) Laplace transform on some neighborhood of $\infty$, and are thus the same \cite[Theorem~8.4]{bhattacharya}. It follows that $\LL_1=\LL_2$. 
%Če je $\LL_1=\LL_2\equiv 0$, sledi $\LL_1(\{\infty\})=1=\LL_2(\{\infty\})$. Sicer je $L_{\LL_1}(0)=L_{\LL_2}(0)>0$, in potem imata zakona $\LL_1(\cdot\cap [0,\infty))/L_{\LL_1}(0)$ ter $\LL_2(\cdot\cap [0,\infty))/L_{\LL_2}(0)$, zožena na Borelove podmnožice $[0,\infty)$, isto končno Laplaceovo transformiranko, in sta torej enaka (glej \cite[Theorem~8.4]{bhattacharya}). Spet sledi $\LL_1=\LL_2$. 
\end{remark}
We now prove Proposition~\ref{proposition:karakterizacija}. Briefly, the idea is to reduce the independence of $R_0$ and $R_1$ to the factorization of their Laplace transforms. This yields a functional equation for the Laplace transforms of the laws of $T_1$ and $T_2$. The latter is in turn analysed using the methods of regular variation, a technique that may be of independent interest. 
\begin{proof}\label{proof}
Let $\fii:=L_{\LL(T_1)}$, respectively $\phi:=L_{\LL(T_2)}$,  be the Laplace transform of the law of $T_1$, respectively $T_2$. 
\label{indeps}
It follows from the relevant independences, i.e. from the fact that the law $\LL((X,T))$ of $(X,T)$ on the space $(\prod_{i\in \mathbb{N}}\{0,1\})\times (\prod_{i\in \mathbb{N}}[0,\infty])$ is given by the product law $\LL((X,T))=(\bigtimes_{i\in \mathbb{N}}\Ber(p))\times (\LL(T_1)\times (\bigtimes_{i\in \mathbb{N}_{\geq 2}}\LL(T_2)))$, from the a.s. equality $\Omega=(\cup_{k\in \mathbb{N}}\{X_1=0\}\cap \{L_1=k+1\})\cup (\cup_{k\in \mathbb{N}}\{X_1=1\}\cap \{L_0=k+1\})$, and from the countable additivity of mathematical expectation, that, for $\{\lambda,\mu\}\subset [0,\infty)$, on the one hand:
$$\EE[e^{-\lambda R_1-\mu R_0}\mathbbm{1}(R_1<\infty,R_0<\infty)]$$ $$=\fii(\lambda+\mu)\left[(1-p)\sum_{k=1}^\infty(1-p)^{k-1}p\phi(\lambda)^k+p\sum_{k=1}^\infty p^{k-1}(1-p)\phi(\mu)^k\right]$$
$$=\fii(\lambda+\mu)p(1-p)\frac{\phi(\lambda)+\phi(\mu)-\phi(\lambda)\phi(\mu)}{(1-(1-p)\phi(\lambda))(1-p\phi(\mu))},$$
and on the other hand: 
$$\EE[e^{-\lambda R_1}\mathbbm{1}(R_1<\infty)]=\fii(\lambda)\left[p+(1-p)\sum_{k=1}^\infty (1-p)^{k-1}p\phi(\lambda)^k\right]=\frac{p\fii(\lambda)}{1-(1-p)\phi(\lambda)},$$
analogously:
$$\EE[e^{-\mu R_0}\mathbbm{1}(R_0<\infty)]=\frac{(1-p)\fii(\mu)}{1-p\phi(\mu)}.$$

The class of bounded functions $\mathcal{K}:=\{e^{-\lambda\cdot}\mathbbm{1}_{[0,\infty)}:\lambda\in [0,\infty)\}$ is closed under multiplication and generates the Borel $\sigma$-field on $[0,\infty]$. The functional monotone class theorem hence implies that the independence of $R_1$ and $R_0$ is equivalent to \footnotesize
\begin{equation}\label{eq:independence-cond}
\EE[e^{-\lambda R_1-\mu R_0}\mathbbm{1}(R_1<\infty,R_0<\infty)]=\EE[e^{-\lambda R_1}\mathbbm{1}(R_1<\infty)]\EE[e^{-\mu R_0}\mathbbm{1}(R_0<\infty)],\quad \{\lambda,\mu\}\subset [0,\infty). 
\end{equation}\normalsize
Indeed, the latter condition is clearly necessary. To see how the functional monotone class theorem intervenes in the proof of the sufficiency, note that for a fixed $\mu\in [0,\infty)$, the class of bounded measurable functions $f:[0,\infty]\to \mathbb{R}$ for which $\EE[f(R_1)e^{-\mu R_0}\mathbbm{1}(R_0<\infty)]=\EE[f(R_1)]\EE[e^{-\mu R_0}\mathbbm{1}(R_0<\infty)]$ is a vector space over $\mathbb{R}$ closed under nondecreasing limits of nonnegative functions. Since it contains the class $\mathcal{K}$ and also $\mathbbm{1}_{[0,\infty]}$, so by monotone class, $\EE[f(R_1)e^{-\mu R_0}\mathbbm{1}(R_0<\infty)]=\EE[f(R_1)]\EE[e^{-\mu R_0}\mathbbm{1}(R_0<\infty)]$ prevails for all bounded measurable $f:[0,\infty]\to\mathbb{R}$. With this having been established, it remains to repeat the preceding argument essentially \emph{verbatim}, except that now with a fixed bounded measurable $f:[0,\infty]\to \mathbb{R}$, and for the class of bounded measurable $g:[0,\infty]\to \mathbb{R}$ for which $\EE[f(R_1)g(R_0)]=\EE[f(R_1)]\EE[g(R_0)]$. 

Using the computations from the beginning of the proof, after some algebraic rearrangement, \eqref{eq:independence-cond} rewrites into
\begin{equation}\label{eq:functional-1}
\fii(\lambda+\mu)\left[\phi(\lambda)+\phi(\mu)-\phi(\lambda)\phi(\mu)\right]=\fii(\lambda)\fii(\mu),\quad \{\lambda,\mu\}\subset [0,\infty).
\end{equation}

The sufficiency of the conditions of the proposition may now be checked as follows. Under \ref{kara:4} $N=0$ a.s., hence $R_0$ and $R_1$ are equal to $\infty$ a.s. and so trivially independent. Under \ref{kara:3}, a.s. $N$ is zero up to $\kappa$ and then jumps to $\infty$ at $\kappa$, whence $R_0$ and $R_1$ are both equal to $\kappa$ a.s., again trivially independent. \ref{kara:2} is the case of a deterministically delayed (by $\kappa$) HPP, in which case the independence of $R_0$ and $R_1$ is well-known, as we have noted.\footnote{Indeed, in all the previous three cases, we see that even $N^0$ and $N^1$ are independent (viz. Remark~\ref{remark}\ref{remark:ii}).} \ref{kara:0}. The deterministic delay by $\kappa$ does not affect independence; we may assume $\kappa=0$. Then $\fii\equiv 1-q_0^2$, $\phi\equiv1-q_0$ and \eqref{eq:functional-1} becomes $(1-q_0^2)[2-2q_0-(1-q_0)^2]=(1-q_0^2)^2$, which holds true. \ref{kara:1}. Again without loss of generality $\kappa$ is set equal to $0$; similarly the scaling of time by the factor $\alpha$ is immaterial to independence, and we may assume $\alpha=1$. In that case we identify $\fii(\lambda)=\frac{1-q_0^2}{1-q_0^2e^{-\lambda}}$ and $\phi(\lambda)=1-q_0+q_0e^{-\lambda}\frac{1-q_0^2}{1-q_0^2e^{-\lambda}}=\frac{(1-q_0)(1+q_0e^{-\lambda})}{1-q_0^2e^{-\lambda}}$, $\lambda\in [0,\infty)$. Tedious but straightforward algebraic manipulations then yield \eqref{eq:functional-1}.

We now prove necessity of the conditions. Set $q_0:=\PP(T_2>0)$. 

Assume $\PP(T_1<\infty)>0$ (otherwise we get \ref{kara:4}) and hence $\fii>0$. We see from \eqref{eq:functional-1} that then  $\phi\not\equiv 0$, hence $\PP(T_2<\infty)>0$, equivalently $\phi>0$. %; in nato (2.) (pri $\lambda=0$) da je $\PP(T_1<\infty)=1$, če je $\PP(T_2<\infty)=1$. %ter končno (3.) pri $\lambda=\mu=0$, da je $\PP(T_1<\infty)=1$, v posebenem $\phi>0$. 
Also from \eqref{eq:functional-1}, for each $\mu\in [0,\infty)$, there exists the limit $$\lim_{\lambda\to\infty}\frac{\fii(\lambda+\mu)}{\fii(\lambda)}=\frac{\fii(\mu)}{1-q_0+\phi(\mu)q_0}\in (0,1].$$ From the characterization of regular variation  \cite[Theorem~1.4.1]{bgt} for the function $\fii\circ \ln\vert_{[1,\infty)}$ it follows that there exists a $\rho \in \mathbb{R}$, for which $$\frac{\fii(\ln r)}{1-q_0+\phi(\ln r)q_0}=r^\rho\text{ for all }r \in [1,\infty);$$ necessarily $\kappa:=-\rho \in [0,\infty)$. In other words $$\fii(\mu)=e^{-\kappa \mu}(1-q_0+\phi(\mu)q_0)\text{ for }\mu\in [0,\infty). $$

If $\PP(T_2=\infty)=q_0$, equivalently if $\LL(T_2)=(1-q_0)\delta_0+q_0\delta_\infty$, then $\phi\equiv 1-q_0$, and we obtain $\fii=e^{-\kappa \cdot}(1-q_0^2)$, whence from the injectivity of the Laplace transform, $\LL(T_1)=(1-q_0^2)\delta_\kappa+q_0^2\delta_\infty$ -- that is to say, we obtain \ref{kara:3} and \ref{kara:0}, according as to whether $q_0=0$ or $q_0>0$.

Assume now $\PP(T_2=\infty)\ne q_0$, equivalently  $\PP(T_2=\infty)<q_0$, in particular $q_0>0$ and $\phi>1-q_0$. 
%Če je $q_0=0$, ekvivalentno $\PP(T_2=0)=1$, je iz $\fii=e^{-\kappa \cdot}(1-q_0+\phi q_0)$ ter iz injektivnosti Laplaceove transformacije, $\PP(T_1=\kappa)=1$. Naj bo sedaj $q_0>0$. Če je $\PP(T_2=\infty)=q_0$, ekvivalentno $\LL(T_2)=(1-q_0)\delta_0+q_0\delta_\infty$, je $\phi\equiv 1-q_0$, in dobimo $\fii=e^{-\kappa \cdot}(1-q_0^2)$. Mera $\LL(T_1)/(1-q_0^2)$, od koder $\LL(T_1)=(1-q_0^2)\delta_\kappa+q_0^2\delta_\infty$. 

The functional equation \eqref{eq:functional-1} may then be rewritten in the form \footnotesize
\begin{equation*}
(1-q_0+q_0\phi(\lambda+\mu))\left[\phi(\lambda)+\phi(\mu)-\phi(\lambda)\phi(\mu)\right]=(1-q_0+q_0\phi(\lambda))(1-q_0+q_0\phi(\mu)),\quad \{\lambda,\mu\}\subset [0,\infty).
\end{equation*}\normalsize
Introducing the substitution $\xi:=(\phi-(1-q_0))/q_0$, we obtain

\begin{equation*}
\xi(\lambda)\xi(\mu)=\xi(\lambda+\mu)(1-q_0^2+q_0^2(\xi(\lambda)+\xi(\mu)-\xi(\lambda)\xi(\mu))),\quad \{\lambda,\mu\}\subset [0,\infty).
\end{equation*}
Another substitution $\psi:=\xi^{-1}-1$ yields 
\begin{equation}\label{eq:functional-2}
\psi(\lambda+\mu)=\psi(\lambda)+\psi(\mu)+\psi(\lambda)\psi(\mu)(1-q_0^2),\quad \{\lambda,\mu\}\subset [0,\infty).
\end{equation}
 
When $q_0=1$, this is  Cauchy's functional equation for $\psi$. The latter being a monotone function, it follows that there exists an $\alpha\in \mathbb{R}$, such that $\psi(\mu)=\alpha\mu$ for all $\mu\in [0,\infty)$. From $\PP(T_2>0)=q_0>0$ we have of course $\alpha\ne 0$ and $\theta:=\alpha^{-1}\in (0,\infty)$. So in this case, for all $\mu\in [0,\infty)$, $\phi(\mu)=\xi(\mu)=(1+\psi(\mu))^{-1}=\theta/(\theta+\mu)$. %V posebnem je $\phi(0)=1$ in zato $ \PP(T_2<\infty)=1$. 
Then the injectivity of the Laplace transform implies that $T_2\sim\Exp(\theta)$. Similarly it follows that $T_1-\kappa\sim \Exp(\theta)$. In other words, we have case \ref{kara:2}.

Finally we are left with the case of $q_0<1$.  When so, we get from \eqref{eq:functional-2} and $\lim_\infty\psi=\infty$ that  for each $\mu\in [0,\infty)$ there exists the limit  $$\lim_{\lambda\to\infty}\frac{\psi(\lambda+\mu)}{\psi(\lambda)}=1+\psi(\mu)(1-q_0^2)\in [1,\infty).$$ The characterization of regular variation gives the existence of an $\alpha\in \mathbb{R}$, such that $$1+\psi(\mu)(1-q_0^2)=e^{\alpha\mu}\text{ for }\mu\in [0,\infty);$$ necessarily $\alpha\in (0,\infty)$. In other words we obtain  $$\phi(\mu)=1-q_0+q_0\frac{1-q_0^2}{e^{\alpha \mu}-q_0^2}\text{ and }\fii(\mu)=e^{-\mu\kappa}\frac{1-q_0^2}{1-q_0^2e^{-\alpha \mu}}\text{ for }\mu\in [0,\infty).$$ %V posebnem je $\PP(T_2<\infty)=\phi(0)=1$. 
The injectivity of the Laplace transform finally implies that $\LL(T_2/\alpha)=(1-q_0)\delta_0+q_0\geom_\mathbb{N}(1-q_0^2)$ and $\LL((T_1-\kappa)/\alpha)=\geom_{\mathbb{N}_0}(1-q_0^2)$, viz. case \ref{kara:1}.
\end{proof} 
As mentioned, Proposition~\ref{proposition:karakterizacija} has as its corollary Theorem~\ref{theorem}.

\noindent \emph{Proof of Theorem~\ref{theorem}.} The conditions that $0$ belong to the support of the law of $T_1$, and that $T_2$ be non-arithmetic or else the law of $T_1$ have no atom in zero, exclude the cases \ref{kara:4}-\ref{kara:3}-\ref{kara:0}-\ref{kara:1} and force $\kappa=0$ in \ref{kara:2} of Proposition~\ref{proposition:karakterizacija}. Clearly the same transpires also when $N$ is ordinary.  \qed\label{proof:thm}

\bibliographystyle{amsplain}
\bibliography{Biblio_HPP-characterization}
\end{document}